\documentclass{article}

\usepackage[dvipsnames]{xcolor}
\usepackage{amsmath, amsthm, amssymb, amsfonts}
\usepackage{thmtools}
\usepackage{graphicx}
\usepackage{setspace}
\usepackage{appendix}
\usepackage{geometry}
\usepackage{float}
\usepackage{tikz-cd}
\usepackage{hyperref}
\usepackage[utf8]{inputenc}
\usepackage[english]{babel}
\usepackage{framed}
\usepackage{relsize}
\usepackage{tcolorbox}
\usepackage{indentfirst}
\usepackage{dynkin-diagrams}
\usepackage{tikz}
\usepackage{tikzsymbols}
\usepackage{biblatex}

\newcommand\R{\ensuremath{\mathbb{R}}}
\newcommand\C{\ensuremath{\mathbb{C}}}

\newcommand\mf{\mathfrak}
\newcommand\Z{\ensuremath{\mathbb{Z}}}

\newcommand\g{\mathfrak{g}}

\newcommand\Hom{\text{Hom}}
\newcommand\h{\mathfrak{h}}
\newcommand\st{\text{ }|\text{ }}
\newcommand{\DynkinBn}[1]

\title{\textbf{Soergel bimodules and Kazhdan-Lusztig polynomials}}
\author{ ethan eugene wynner}
\date{University of California, Berkeley}

\begin{document}

\maketitle

\section*{Abstract}

This paper presents an exposition of Soergel bimodules, the categorification of Hecke algebras, and their applications to various topics in Kazhdan-Lusztig theory. We ultimately exposit a few of Soergel's main results, which allowed him to give alternative proofs, using his theory, of the Kazhdan-Lusztig conjectures. This paper should be viewed as a (very) condensed exposition following the work of primarily Elias, Makisumi, Thiel, and Williamson in their lovely book \textit{Introduction to Soergel Bimodules}\ref{ref2}.  It is also important to note that much of Soergel's theory has been omitted from this paper in light of brevity, as well as the depth of proofs and lead-up information. Many proofs concerning covered material are left out of this paper, again, due to their length, as well as their residence in many of the texts listed in the references section. \ref{references}

\newpage
\tableofcontents

\newpage
\begin{center}
    \section{Coxeter Systems, Hecke algebras, and Kazhdan-Lusztig preliminaries}
\end{center}
\subsection{Coxeter Systems}

Suppose $W$ is a group, and $S \subset W$ is a finite subset consisting of generators of $W.$ Let $(m_{st})_{s,t \in S}$ be a matrix satisfying 
$$m_{st} = m_{ts} \in \{1,2,\ldots\} \cup \{\infty\} $$ with $m_{ss} = 1$ for every $s \in S$. We say that the pair $(W,S)$ is a \textit{Coxeter system}, and $W$ a \textit{Coxeter group}, if $W$ admits the presentation 
$$W = \langle s \in S \st (st)^{m_{st}} = e \text{ for any } s,t \in S, m_{st} < \infty\rangle,$$
where $e$ denotes the identity element in $W$. It turns out that the value $m_{st}$ is the order of the element $st \in W$, and, if $m_{st}$ is not finite, no correlation exists between $s$ and $t$. One will often see generators $s \in S$ submitted to the condition that it is involutive, i.e. $s^2 = e$. These generators are sometimes called reflections, and so it may not be a surprise that 
 \subsubsection{Remark:}  
  A Weyl group $W$ is a Coxeter group. \\
  Indeed, if $\g$ is a complex semisimple Lie algebra, and $\h$ its Cartan (maximal toral) subalgebra, $W$ is generated by reflections across root hyperplanes, given by operators $s_\alpha: \h^* \rightarrow \h^*$, with
  $$s_\alpha(\lambda) = \lambda - \dfrac{2(\lambda, \alpha)}{(\alpha,\alpha)}\alpha.$$ Here $\alpha \in R$, where $R$ denotes the root system of $\g$, and ( , ) is the euclidean inner product.\ref{app1}

  \subsection{Hecke algebras}
  Fix a Coxeter system $(W, S),$ and an indeterminate parameter $v$ over $\Z.$ Let $\Z[v]$ be the polynomial ring generated by $v,$ and let $\Z[v, v^{-1}]$ be its localization. The \textit{Hecke algebra} $\mathcal{H}$ of $W$ is a unital associative  $\Z[v, v^{-1}]$-algebra generated by the set $\{\delta_s \st s \in S\}$ under the quadratic relation 
  $$\delta^2_s = (v^{-1} - v)\delta_s+1,$$
  and the braid relation 
  $$\underbrace{\delta_s\delta_t\delta_s\ldots}_{m_{st}} = \underbrace{\delta_t\delta_s\delta_t\ldots}_{m_{st}} $$
  for all simple reflections $s,t \in S.$ If $R$ is any commutative ring, the \textit{specialization} of $\mathcal{H}$ is a ring homomorphism 
  $$\varphi: \Z[v, v^{-1}] \rightarrow R, \text{    }v \mapsto \varphi(v)$$
  along with the space
  $$\mathcal{H}_\varphi = R \otimes_{\Z[v, v^{-1}]} \mathcal{H}, $$
  uniquely determined by $R$ and the mapping $v \mapsto \varphi(v).$

  \subsubsection{Example}
 Let $R = \Z,$ with specialization $ \varphi: v \mapsto 1.$ Then $\mathcal{H}_\varphi = \Z \otimes_{\Z[v, v^{-1}]} \mathcal{H}$ is isomorphic to the group algebra $\Z[W]$ with $\delta_s \mapsto s, \text{  }  \forall s \in S.$ 
 \\ \textit{Proof:} Define 
 $$\varphi: \Z[v, v^{-1}] \rightarrow \Z, \text{    }v \mapsto 1.$$
  Then 
 $$\delta^2_s = (1^{-1} - 1)\delta_s + 1 = 1 = \delta_{id}.$$
 Hence $\delta^2_s$ satisfies the property $s^2 = e$ in $W,$ and the braid relation holds by definition. We thus establish the bijection 
 $$\{\delta_w \st w \in W\} \longleftrightarrow \{w \in W\},$$
 from which the bijection 
 $$\Z \otimes_{\Z[v, v^{-1}]} \mathcal{H} \longleftrightarrow \Z[W]$$
arises. Since this mapping preserves the defining relations in  $\Z[W],$ it is an isomorphism.  \Laughey

 \subsubsection*{}
 
 In general, under the specialization $v \mapsto 1,$ the Hecke algebra recovers the group algebra $\Z[W],$ hence we consider $\mathcal{H}$ to be a \textit{deformation} of the $\Z[W].$

 \subsection{Kazhdan-Lusztig polynomials, and a basis for the Hecke algebra}

 Let $x \in W$, and define 
 $$\delta_x = \delta_{s_1}\delta_{s_2}\cdots\delta_{s_m},$$
 for a reduced expression $\Tilde{x} = (s_1, s_2, \ldots s_m).$
 We call the involution $\mathcal{H} \rightarrow \mathcal{H},$ $h \mapsto \bar{h}$ defined by 
 $$\bar{\delta} = \delta^{-1}_s = \delta_s + (v - v^{-1}),$$
 the \textit{Kazhdan-Lusztig involution.} We define the \textit{Kazhdan-Lusztig basis} of $\mathcal{H}$ to be the set $\{b_x \st x \in W\},$ where elements $b_x$ are self-dual and of the (degree-bounded) form 
 $$b_x = \delta_x + \sum_{y < x}h_{y,x}\delta_y,$$
 where $y < x$ is in the sense of Bruhat decomposition, and $h_{y,x} \in v\Z[v]$. These coefficients $h_{y,x}$ are called \textit{Kazhdan-Lusztig polynomials.} 

 \subsubsection{Lemma:} 
 \textit{A Kazhdan-Lusztig basis is unique.}
 \\ A tidy proof of this lemma can be found in \ref{ref2} 

 \subsubsection{Lemma:} \textit{For every $s \in S,$ the Kazhdan-Lusztig basis element $b_s$ has the form $b_s = \delta_s + v$}
 \\ \textit{Proof:} $\delta_s + v$ satisfies the degree bounded form by definition, and we have 
 $$\overline{\delta_s + v} = \delta^{-1}_s + v^{-1} = \delta_s + (v - v^{-1}) + v^{-1} = \delta_s + v,$$
 Uniqueness of the Kazhdan-Lusztig basis gives us the desired result. \Laughey

\begin{center} \section{Soergel Bimodules} \end{center}

\subsection{Easing into it}

Let $(W,S)$ be a Coxeter system, and let $V$ be its geometric representation\ref{app2} over $\R$. Let $\varPi = \{\alpha_s \st s \in S\}$ be the set of simple roots, forming a basis of $V$ indexed by simple reflections. Let 
$$R = \text{Sym}(V) = \bigoplus_{i \in \Z_{\ge 0}}\text{Sym}^i(V)$$
be the symmetric algebra of $V,$ understood here as a $\Z$-graded algebra where deg($V$) = 2. Then $R$ can be stated equivalently as $\R[\varPi]$ with grading deg($\alpha_s)= 2$.

\subsubsection{Theorem (Chevalley-Shephard-Todd)}
Let $I \subset S$ be a subset such that $W_I = \langle I\rangle \subset W$ is finite. Then the ring of  $W_I$-invariants 
$$R^{W_I} = \{x \in R \st w \cdot x = x, \forall w \in W\}$$
is a polynomial ring, and $R$ is a $\Z$-graded free module over $R^{W_I}$ of finite rank.

\subsubsection{Example}
Let $W$ be the Weyl group of $\g = \mf{sl}_2,$ i.e. $W = S_2,$ such that $S = \{s\},$ i.e. $W_I$ is generated by a single reflection. Letting $R = \R[\alpha]$, the $W$-action is given by $s(\alpha) = -\alpha.$ Thus $R^s = \R[\alpha^2].$ It follows that $R$ has the decomposition
$$R = R^s \oplus (R^s \cdot \alpha).$$
More generally, in a Coxeter system any polynomial $P \in R$ can be written $P = a+ b\alpha_s,$ where $a,b \in R^s.$ Since $\deg(\alpha_s) = 2,$ we have $R \cong R^s \oplus R^s(-2)$ as graded bimodules. 

\subsection{Bott-Samelson bimodules}
Let $B_s$ for $s \in S$ denote the graded $R$-bimodule
$$B_s = R \otimes_{R^s}R(1),$$
where $R(1)$ is obtained by shifting down by 1. That is, $M(1)^i \equiv M^{i-1}.$ Elements in $B_s$ are expressable as sums 
$$\sum_i f_i \otimes g_i, \text{ }f_i, g_i \in R.$$
Given a reduced expression $\underbar{$w$} = (s_1, \ldots s_n).$ We define the Bott-Samelson bimodule to be the graded $R$-bimodule 
$$BS(\underbar{$w$}) = B_{s_1}B_{s_2}\cdots B_{s_n},$$
equivalently, 
$$BS(\underbar{$w$}) = R \otimes_{R^{S_1}}R \otimes_{R^{S_2}} \cdots \otimes_{R^{S_n}} R(l(w)).$$
Thus, we express our elements of $BS(\underbar{$w$})$ as 
$$\sum_i f_i \otimes g_i \otimes \cdots \otimes h_i,$$
with of course, $f_i, g_i, \ldots h_i \in R.$ We can now define the titular bimodule of this paper: A \textit{Soergel bimodule} is a direct summand of a finite direct sum of grading shifts of Bott-Samelson bimodules. The category of Soergel bimodules $\mathbb{S}$Bim is the strictly full smallest subcategory of $R$-gbim\ref{app3} (the category of graded $R$-bimodules) containing both $R$ and $B_s$ for every $s \in S$, and is closed under tensor products, direct sums/summands, and grade shifting. 

\subsubsection{Example of Soergel bimodule}
Consider the polynomial ring $ R = \mathbf{k}[x,y,z]$ where $\mathbf{k}$ is a field not of characteristic $2,$ and consider the action of $S_3$ on $R,$ characterized by simple reflections $s$ and $r,$ where 
$$s \cdot \mu(x,y,z) = \mu(y,x,z)$$
and 
$$r \cdot \mu(x,y,z) = \mu(x,z,y).$$
Hence $R^s$ is the polynomial ring $\mathbf{k}[x+y,xy,z],$ and $R^r$ is $\mathbf{k}[x+y+z,xy+xz+yz, xyz].$ We take $\deg(x,y,z) = 2.$ Then $R(1)$ is $R$ shifted down by 1, which we obtain by declaring the degree of each polynomial reduced by one. Specifically, if $B = \bigoplus_{i \in \Z}B^i$ is a $\Z$-graded object, we declare 
$$B(m) = \bigoplus_{i \in \Z}B(m)^i = \bigoplus_{i \in \Z}B^{m+i}.$$
Then the $\Z$-graded $R$-bimodule $R$ is a Soergel bimodule\ref{ref9}. 

\subsubsection*{}
Have you ever wondered, "What are the morphisms like in the category of Soergel bimodules?"? To satisfy this musing, we take for granted that any Soergel bimodule is graded free as a left or right $R$-module, and the assumption that morphisms between Soergel bimodules are homogenous with degree 0. Formally, we describe 
$$\Hom_{\mathbb{S}\text{Bim}}(B, B') = \Hom_{R\text{-gbim}}(B, B') $$ where 
$$\Hom_{\mathbb{S}\text{Bim}}(B, B') \simeq \Hom_{\mathbb{S}\text{Bim}}(B, B'(k)).$$
Notice how this juxtaposes Bott-Samelson bimodules, which are not closed under grading shifts, as $\mathbb{S}\text{Bim}$ is. To do some housekeeping, we introduce the category of Bott-Samelson bimodules, $\mathbb{B}\mathbb{S}\text{Bim}$, as a monoidal category which is not closed under grading shifts. Hence morphisms $B \rightarrow B'$ in this category are characterized 
$$\Hom_{\mathbb{B}\mathbb{S}\text{Bim}}(B,B') = \bigoplus_{k \in \Z} \Hom_{R\text{-gbim}}(B,B'(k))$$
as a graded vector space, where the RHS summands are the spaces of homogenous bimodule maps of degree k. 

Recall what it means to be an \textit{indecomposable} object $M$ in an additive category $\mathcal{A}.$ Specifically, if $M$ is indecomposable, then it cannot be expressed as a direct sum 
$$M = M' \oplus M'',$$
where $M',M'' \ne 0.$ We introduce the following lemma to make a coming example more robust. 

\subsubsection{Lemma}
\textit{Suppose $M$ is a graded $R$-bimodule generated by a homogenous element $m \in M.$ Then $M$ is indecomposable.}
\\ \textit{Proof:} We have $M = RmR.$ Let $d = \deg(m).$ Then $M^d = R^0mR^0,$. Suppose that $M = L \oplus N.$ Then $M^d = L^d \oplus N^d,$ so we can say WLOG $m \in L.$ Then $M = R \cdot m \cdot R \subset L,$ which implies that $N^d = 0,$ so that $M = RmR \subset N,$ i.e. $N = 0.$   \Laughey

\subsubsection{Example}
Now we can formulate another example in type $A_1,$ wherein $W = S_2, $ generated by the simple reflection $s.$ By the previous lemma, $R$ and $B_s$ are indecomposable. Hence 
$$B^2_s = B_sB_s \simeq R \otimes_{R^s} R \otimes_{R^s} R(2) \simeq R \otimes_{R^s} (R^s \otimes R^s(-2)) \otimes_{R^s} R(2)$$
$$\simeq R\otimes_{R^s} R(2) \oplus R \otimes_{R^s} R \simeq B_s(1) \oplus B_s(-1).$$

The reason that this example is interesting is that it corresponds directly to the Hecke algebra $\mathcal{H}$ of type $A_1.$ Indeed, where $b_s$ is the Kazhdan-Lusztig basis element in $\mathcal{H},$ we have 
$$b_sb_s = (v+v^{-1})b_s = b_sv + b_sv^{-1}.$$
Other glorious examples of this correspondence in different types of Coxeter systems can be found in [Elias et al]
\subsubsection*{}
This gives rise to the following

\subsubsection{Remark}
Whenever $s \ne t$ and $m_{st} \ne \infty,$ $R$ is generated by subrings $R^s$ and $R^t.$ We observe that 
$$B_sB_t \simeq R \otimes_{R^s} R \otimes_{R^t} R(2) $$
and 
$$B_tB_s \simeq R \otimes_{R^t} R \otimes_{R^s} R(2)$$
are indecomposable, since they are not generated by the 1-tensor $1 \otimes 1 \otimes 1,$ which in turn implies that $B_s$ and $B_t$ are not isomorphic. 
\\ \subsubsection*{Subremark}
 It is important now to collect ourselves and remember that $R = \text{Sym}(V),$ where $V$ is the geometric representation of $W.$ As a polynomial ring, $R$ is generated by its linear terms, i.e. the vectors in $V,$ whereas the linear terms in $R^s, R^t,$ are the vectors in $V$ that are fixed by $s,t$ respectively. Rather, the reflecting hyperplanes $H_s$ and $H_t$ are perpindicular to $\alpha_s$ where $s\ne t,$ and as such, they together span $V.$ 

 \subsection{Grothendieck groups,  Soergel's categorization}
 Probably the most important technical attribute of Soergel bimodules is their role as the "categorification" of the Hecke algebra, and their connection to the Kazhdan-Lusztig basis of said Hecke algebra. In pursuit of properly expositing this, we will interface now with one of the many beautiful constructions in math, the Grothendieck group. 
\\ The \textit{split Grothendieck group} $[\mathbb{S}\text{Bim}]_\oplus$ of $\mathbb{S}\text{Bim}$ is an abelian group generated by classes $[B]$ of Soergel bimodules in $\mathbb{S}\text{Bim}.$ We equip the following simple yet universally powerful condition that 
$$[B] = [B'] + [B'']$$
where $B \simeq B' \oplus B''.$ Since $\mathbb{S}\text{Bim}$ is a monoidal category, we can give $[\mathbb{S}\text{Bim}]_\oplus$ the structure of a ring by setting $[B][B'] = [BB'].$ Moreover, by setting $v[B] = [B(1)],$ we give $[\mathbb{S}\text{Bim}]_\oplus$ the structure of a $\Z[v^{\pm 1}]$-algebra. 

\subsubsection{Definition}

The \textit{character} $\chi$ of a Soergel bimodule $B$ is defined as the element 
$$\chi(B) = \sum_{x \in W}v^{l(x)}h_x(B)\delta_x$$
of $\mathcal{H}.$ We have $h_{id}(B_s) = v^{+1}$ and $h_s(B_s) = v^{-1},$ and thus $\chi(B_s) = v\delta_{id} + v - v^{-1}\delta_s = v + \delta_s.$ Thus, wonderfully, $\chi(B_s) = b_s $ for every $s \in S.$ It follows that $\chi(B \oplus B') = \chi(B) + \chi(B')$ and $\chi(B(1)) = v\chi(B).$ Since $v[B] = [B(1)]$ in the Grothendieck group, $\chi$ induces a $\Z$-linear, but moreover $\Z[v, v^{-1}]$-linear map 
$$\chi: [\mathbb{S}\text{Bim}]_\oplus \rightarrow \mathcal{H}.$$
We can now introduce these very important results : 

\subsubsection{Soergel's Categorification theorem}
\begin{enumerate}
    \item \textit{There exists a $\Z[v^{\pm 1}]$-algebra homomorphism}
    $$\phi: \mathcal{H} \rightarrow [\mathbb{S}\text{Bim}]_\oplus, \text{       }b_s \mapsto [B_s]$$
    \textit{for each $s \in S$.}

    \item \textit{$W$ is in bijection with the set of indecomposable objects in $\mathbb{S}\text{Bim}$ mod shift and isomorphism. Formally, }
    $$W \longleftrightarrow \{\text{indecomposable objects} \in \mathbb{S}\text{Bim}\}/\simeq, (k), \text{     } w \mapsto B_w, $$
    \textit{where $B_w$ are direct summands in the Bott-Samelson bimodules $BS(\underline{w})$}.

    \item \textit{$\chi$ descends to a $\Z[v^{\pm 1}]$-module homomorphism }
    $$\chi: [\mathbb{S}\text{Bim}]_\oplus \rightarrow \mathcal{H}$$
    \textit{which is inverse to $\phi$. Thus $\phi$ and $\chi$ are isomorphisms, and we have $[\mathbb{S}\text{Bim}]_\oplus \simeq \mathcal{H}.$} 
\end{enumerate}

\subsubsection{Soergel's Hom formula}

\textit{For any Soergel bimodules $B, B'$, the graded $R$-bimodule hom space $\Hom^\bullet_{\mathbb{S}\text{Bim}}(B, B')$ is free as a graded left (resp., right) $R$-module. Moreover, }
$$\text{gr(rank)}\Hom^\bullet_{\mathbb{S}\text{Bim}}(B, B') = (\chi(B), \chi(B')),$$
\textit{where (  ,  ) is the standard form on $\mathcal{H}.$} 

\subsubsection*{}
Soergel's Hom formula motivates what we call \textit{Soergel's conjecture,} which was proved for arbitrary Coxeter systems (Soergel proved for Weyl groups) by [Elias / Williamson]. We state it here: 

\subsubsection{Soergel's conjecture}
\textit{For every $x \in W$, $\chi(B_x) = b_x.$ Equivalently, the Kazhdan-Lusztig polynomial $h_{x,y} = h_x(B_y)$.}

\subsubsection{Important Remark}
A positive result of Soergel's Conjecture has a very important implication in the larger theory. That is, it implies the \textit{Kazhdan-Lusztig positivity conjecture,} which asserts that the coefficients of $h_{x,y}$ are non-negative. Indeed, $h_x(B)$ lies in $\Z_{\ge 0}[v^{\pm 1}]$ for any Soergel bimodule $B.$

\subsubsection*{}
We will now go further into some of the historical context in Lie theory and Category $\mathcal{O},$ with a little algebraic geometry. It is important for the reader to note here that the history/applications of this theory are wide ranging, and greatly exceed the scope of this paper. The interested reader can find more comprehensive texts in \ref{references}.

\newpage
\section{Historical Interlude: The Kazhdan-Lusztig conjectures}

\subsection{Briefly recognizing some Lie theory}
Let $\g$ be a complex semisimple Lie algebra. We have the decomposition 
$$\g = \mf{n}_+ \oplus \h \oplus \mf{n}_-$$
with $\h$ the Cartan subalgebra of $\g$ and $\mf{n}_\pm = \bigoplus_{\alpha \in R_+} \g_\alpha.$ The \textit{Borel subalgebra} $\mf{b}$ is the subalgebra $\mf{b} = \mf{n}_+ \oplus \h.$ Under the relations
$$hv_\lambda = \langle h, \lambda\rangle v_\lambda \text{        }\forall h \in \h$$
$$xv_\lambda = 0 \text{     } \forall x \in \mf{n_+},$$
we obtain a 1-dimensional representation $\C_\lambda$ of $\mf{b},$ and we define the \textit{Verma module} $M_\lambda$ to be 
$$M_\lambda = U(\g) \otimes_{U(\mf{b})} \C_\lambda.$$

\subsubsection{Definition: Category \texorpdfstring{$\mathcal{O}$}{O}}

\textit{Category} $\mathcal{O}$ is the abelian full subcategory of Rep($\g$) with 
$$V = \bigoplus_{\lambda \in \h^*} V[\lambda]$$
where $V[\lambda]$ are weight spaces\footnote{$V[\lambda] := \{v \in V \st hv = \lambda(h)v \text{   } \forall h \in \h\}$}.  Further, we ask that $V$ is finitely generated, and the action of $\mf{n}_+$ on $V$ is locally finite. Notice that the set of Verma modules $\{M_\lambda \st \lambda \in \h^*\}$ is contained in $\mathcal{O}.$ Furthermore, the set of highest weight irreducible $U(\g)$-modules, denoted $\{L_\lambda \st \lambda \in \h^*\},$ gives a complete set up to isomorphism of simple modules in $\mathcal{O}.$ We define $\mathcal{O}_0$ to be the épaisse\footnote{épaisse, meaning "thick", denotes a nontrivial full category which is closed under direct sums.} subcategory of $\mathcal{O}$ generated by $\{L_{(w \cdot 0)} \st w \in W\}.$ 

\subsubsection{Definition: Soergel module}
We have already spent a good chunk of this paper talking about Soergel bimodules without any discussion of a \textit{Soergel module.} This is for good reason, since Soergel modules are obtained via an a priori construction of Soergel bimodules, using some relevant information in this section. We do this now, given our complex semisimple Lie algebra $\g$ and its Weyl group $W,$ which has the structure of a Coxeter system. The Cartan subalgebra $\h \subset \g$ functions as a realized version of this Coxeter system over $\C.$ Consider Soergel bimodules associated with the $W$-action on the symmetric algebra $\text{Sym}(\h).$ A \textit{right Soergel module} is a graded $R$-module of the form 
$$\C \otimes_{\text{Sym}(\h)}B,$$
where $B$ is a Soergel bimodule and the map $\text{Sym}(\h) \rightarrow \C$ sends $h \mapsto 0.$ So, essentially, we obtain right (resp. left) Soergel modules by annihilating the action of higher degree polynomials on the left. 

\subsection{Back to Kazhdan-Lusztig}
Recall the Kazhdan-Lusztig basis $\{b_w \st w \in W\}$ of the Hecke algebra $\mathcal{H},$ given a Coxeter system $(W,S).$ We want to know what this has to do with a Weyl group $W.$ $\mathcal{H}$ also admits a standard basis $\{\delta_w \st w \in W\},$ and the translation matrix from the KL basis to the standard one is upper triangular with respect to Bruhat order. We have
$$b_x = \sum_{y \le x}h_{y,x}(v)\delta_y,$$
where the Kazhdan-Lusztig polynomials are the recursively defined coefficients $h_{y,x}(v).$ In fact, this construction is reliant only on the combinatorics of the Weyl group. We now have a theorem of Bernstein-Gelfand-Gelfand about Verma modules:

\subsubsection{Theorem}
\textit{If $\Hom(M_\mu, M_\eta) = 1,$ then there exists a dominant weight $\lambda \in \h^*$ and $x,y \in W$ such that}
$$\mu = x \cdot \lambda, \eta = y \cdot \lambda,$$
with $x \ge y$ in Bruhat order. This brings us our seminal conjecture, stated as follows.

\subsubsection{The Kazhdan-Lusztig conjecture}

\textit{ For every $x,y \in W,$ the following holds:}
$$[M_{(y\cdot 0)} : L_{(x \cdot 0)}] = h_{y,x}(1).$$

This conjecture is the focal application of Soergel bimodules. In the near future, we will cover this in more depth. Kazhdan and Lusztig also proved the following formula, the 

\subsubsection{Kazhdan-Lusztig Inversion Formula}
Let $w_0 \in W$ be the longest eleemnt. For every $x,y \in W,$ we have the identity
$$\sum_{z \in W}(-1)^{l(y) + l(x)}h_{z,x}(v)h_{zw_0, yw_0}(v) = \begin{cases}
    1 \text{ if } x = y \\
    0 \text{ else}
\end{cases}$$

Along with the previous conjecture, this implies a similar identity in the Grothendieck group of $\mathcal{O}$:

$$[L_{(y\cdot 0)}] = \sum_{y \le x}(-1)^{l(y) + l(x)}h_{xw_0, yw_0}(1)[M_{(x\cdot 0)}].$$

At this point, we might find ourselves asking some questions about what exactly Kazhdan-Lusztig polynomials are. To try shedding some light, we turn to the following theorem of of Polo \ref{ref11}, stating roughly that any polynomial with positive coefficients is a Kazhdan-Lusztig polynomial in some symmetric group.

\subsubsection{Theorem (Polo)}
\textit{For any monic polynomial $q \in \Z_{\ge 0}[v^2]$, there exist numbers $m$ and  $N$ such that}
$$v_mq = h_{y,x}(v)$$
\textit{for $x,y \in S_N$.}

\subsection{More on positivity}

Recall the conjecture from earlier of Kazhdan-Lusztig, stating that Kazhdan-Lusztig polynomials have non-negative coefficients. The original proof (done by Kazhdan and Lusztig) required the geometric notion of flag varieties. We'll briefly introduce this here.  

\subsubsection{Flag and Schubert varieties}
Recall the decomposition of $\g$ by 
$$\g = \mf{n}_+ \oplus \h \oplus \mf{n}_-$$
Let $G$ denote the simply connected Lie group corresponding to $\g,$ and let $B, U$ be the subgroups of $G$ corresponding to the Borel subalgebra $\mf{b}$ and $\mf{n}_+,$ respectively. Then the \textit{flag variety} is defined to be the quotient $G/B$. Considering $U$-orbits (equivalently, left $B$-orbits), we obtain the \textit{Bruhat decomposition}
$$G/B = \coprod_{w \in W}BwB/B.$$
The \textit{Schubert variety} $X_w$ is the closure of $BwB/B.$ Coinciding with the Bruhat order is the relation $y \le w$ if $X_y \subset X_w.$

\section{Soergel's Proof of the Kazhdan-Lusztig conjecture}

\subsection{The \texorpdfstring{$\mathbb{V}$}{V}-functor}
Soergel's proof makes use of two important functors. One of them, denoted $\mathbb{V} $, connects the subcategory $\mathcal{O}_0$ to Soergel modules, which we have just defined in the previous section. In pursuit of a more mathematical definition, we consider the fact that a $W$-invariant polynomial $f \in R^{W_I}$ passes freely through every tensor product ($\otimes_{R^s})$ in a Bott-Samelson bimodule $BS(\underline{w}).$ This implies two things. One is that $f \in R^{W_I}$ will act the same on the right and left of any Soergel bimodule. Secondly, if $f$ is of strictly positive degree, it will act as zero on any Soergel \textit{module}. A definition here clarifies this acknowledgement: 

\subsubsection{Definition}
If $W$ is finite, let $R^{W_I}_+$ be the graded subspace of $R^{W_I}$ containing everything of strictly positive degree. Let $I_W$ be the homogenous ideal generated by $R^{W_I}_+$. The \textit{coinvariant algebra} of $W$ is the graded algebra 
$$C = R/I_W.$$
Following from above, $I_W$ acts as zero on any Soergel module, hence Soergel modules in general can be viewed as graded right $C$-modules. 

\newpage
Now we turn back to $\mathcal{O}.$ Each $L_\lambda$ has an indecomposable projective cover, which we will denote by $\text{Pr}_\lambda.$ Define $\mathbb{V}$ by $\mathbb{V} := \Hom(\text{Pr}_{(w_0 \cdot 0)}, -).$ Soergel proved the following theorem about $\mathbb{V}:$

\subsubsection{Theorem}
\textit{There is a canonical isomorphism}
$$\text{End}(\text{Pr}_{(w_0\cdot 0)}) \simeq \mathbb{V}(\text{Pr}_{(w_0\cdot 0)}) \simeq C.$$
Using this, we view $\mathbb{V}$ as a functor 
$$\mathbb{V}: \mathcal{O}_0 \rightarrow \text{mod-}C$$
where mod-$C$ are ungraded right $C$-modules. Consequently, 

\subsubsection{Theorem}
\textit{Let} Proj$\mathcal{O}_0 \subset \mathcal{O}_0$ \textit{be the full subcategory of projective modules. Then}
\begin{enumerate}
    \item $\mathbb{V}|_{\text{Proj}\mathcal{O}_0}$ \textit{is fully faithful}.
    \item $\mathbb{V}(\text{Pr}_{(x\cdot 0)}) \simeq \overline{B_x}$ \textit{as ungraded modules for every $x \in W.$}
\end{enumerate}

 \subsection{Translation functors}
Fixing the decomposition of $\g = \mf{n}_+ \oplus \h \oplus \mf{n}_- $, the weight lattice $P \subset \h^*,$ root system\footnote{it is important to clarify that what we mean by $R_\pm$ will always be a root system. Any use of $R$ without a plus or minus subscript will be the symmetric algebra associated to our construction of Soergel (bi)modules, unless explicitly stated otherwise.} $R_\pm$, simple roots $\Pi$, and the Weyl group $W,$ we consider the $W$-action on $\h^*$, given by 
$$w\cdot \lambda = w(\lambda+ \rho) - \rho,$$
for $w \in W, \lambda \in h^*, \rho = \dfrac{1}{2}\sum_{\alpha \in R_+}\alpha.$ Then $\mathcal{O}$ has the decomposition
$$\mathcal{O} = \bigoplus_{\lambda \in \h^* / W}\mathcal{O}_\lambda,$$
with simple objects $\{L_{(w \cdot \lambda)} \st w \in W\}$ in each $\mathcal{O}_\lambda.$ 

\subsubsection{}
Let $\lambda \in \h^*$. Define 
$$\iota_\lambda: \mathcal{O}_\lambda \hookrightarrow \mathcal{O}
, \text{           }\pi_\lambda: \mathcal{O} \twoheadrightarrow \mathcal{O}_\lambda$$
as the inclusion and projection, respectively, associated with $ \mathcal{O}_\lambda.$ 
Let $\lambda, \mu \in \h^*$ be compatible ($\lambda - \mu$ is integral). Let $\eta$ be the unique dominant integral weight in $W_{(\lambda - \mu)}$ Let $V_\eta$ be the finite dimensional irreducible representation of highest weight $\eta.$ The \textit{translation functor} $T^\mu_\lambda$ is the exact (by composition) functor 
$$T^\mu_\lambda := \pi_\mu \circ (V \otimes (-)) \circ \iota_\lambda : \mathcal{O}_\lambda \rightarrow \mathcal{O}_\mu. $$

\subsubsection{Lemma}
\textit{Suppose that $\lambda, \mu \in \h^*$ are compatible. Then the translation functors}
$$\begin{tikzcd}
    \mathcal{O}_\lambda 
    \arrow[r, bend left, "T^\mu_\lambda"] 
    & 
    \mathcal{O}_\mu 
    \arrow[l, bend left, "T^\lambda_\mu"]
\end{tikzcd}$$
\textit{are biadjoint.}
\\\\ \textit{Proof:} If $V_\eta$ is defined as it is above, then $V^*_\eta = V_{(-w_0\eta)}$ where $w_0$ is the longest element in $W.$ Then $T^\lambda_\mu = \pi_\lambda \circ (V^* \otimes (-)) \circ \iota_\mu.$ The desired result follows from the fact that $V \otimes (-)$ and $V^* \otimes (-)$ are biadjoint. \Laughey

\subsubsection{Lemma}
\textit{Translation functors send projective modules to projective modules.}
\\\\ Proof of this fact follows from the previous lemma and a more general result about projective objects in abelian categories, the proof of which has been relegated to the appendix. \ref{app4}

\subsubsection*{}
With our earlier setup for $\g$,  let $E$ be the real span $E = \R \otimes_\Z P \subset \h^*$, equipped by the Killing form with the inner product ( , ). For any root $\alpha \in R_\pm$, define the hyperplane $H_\alpha$, orthogonal to $\alpha$, by 
$$H_\alpha := \{\lambda \in E \st (\lambda + \rho, \alpha) = 0\},$$
or equivalently, $H_\alpha := \{\lambda \in E \st (\lambda + \rho)(\alpha^\vee) = 0\}.$ In $\mathcal{O}_0,$ we will give an important construction. \\\\
Let $s \in S,$ and choose an integral weight $\mu$ such that the dot-action stabilizer of $\mu$ is $\{\text{id}, s\},$ and choose a dominant integral weight $\nu$ such that $\nu - \mu$ is also dominant integral. Define the composition of functors 
$$\Theta_s:= T^0_\nu \circ T^\nu_\mu \circ T^\mu_\nu \circ T^\nu_0.$$
An important property of this construction is the following, wherein we let the Verma module $M_{(w \cdot 0)}$ be denoted by $M_w$.

\subsubsection{Lemma}
\textit{For any $w \in W$, $s\in S$, the following sequence is nonsplitting and short exact:}
$$\begin{tikzcd}
    0 \arrow[r] & M_w \arrow[r] & \Theta_sM_w \arrow[r] & M_{ws} \arrow[r] & 0
\end{tikzcd} \text{ if } ws > w$$
$$\begin{tikzcd}
    0 \arrow[r] & M_{ws} \arrow[r] & \Theta_sM_w \arrow[r] & M_{w} \arrow[r] & 0
\end{tikzcd} \text{ if } ws < w$$

Noticing that $[\Theta_sM_w] = [M_w] + [M_{ws}]$ in the Grothendieck group, we observe the following

\subsubsection{Corollary}
\textit{The group isomorphism}
$$[\mathcal{O}_0] \rightarrow \Z[W], \text{     } [M_w] \mapsto w$$
\textit{intertwines the map induced by $\Theta_s$ with right multiplication by $1 + s$. }

\subsection{Projectives}
Let $\text{Pr}_{(w \cdot 0)}$ be denoted as $\text{Pr}_w$ going forward. 

\subsubsection{Proposition}
\textit{The set} $\{\text{Pr}_w \st w \in W\}$ \textit{is the set of non-isomorphic indecomposable projective objects in} $\mathcal{O}_0.$ \textit{Moreover, every projective object $Q \in \mathcal{O}_0$ is isomorphic to a finite direct sum }
$$Q \simeq \bigoplus_{w \in W}\text{Pr}^{\oplus m(\text{Pr}_w, Q)}_w$$
\textit{for unique multiplicites } $m(\text{Pr}_w, Q) \ge 0$, with $m(\text{Pr}_w, Q) = \dim \Hom(Q, L_w).$ \\\\
\textit{Proof:} Let $Q \in \mathcal{O}_0$ be a projective object. Choose a filtration 
$$0 = N_0 \subset N_1 \subset \ldots \subset N_n = Q$$
with simple successive subquotients. Then we have the natural projection $Q \rightarrow N_n/N_{n-1} \simeq L_x$ for some $x \in W.$ Since $\text{Pr}_x$ is projective, there exist maps $f: \text{Pr}_x \rightarrow Q$ and $g: Q \rightarrow \text{Pr}_x$ such that $gf(\text{Pr}_x) \subset \text{Pr}_x$ is a submodule of $\text{Pr}_x.$ Since $\text{Pr}_x$ is indecomposable, we have $gf(\text{Pr}_x) = \text{Pr}_x.$ Hence $\text{Pr}_x$ is a direct summand of $Q.$ By induction on the length of filtration, we ascertain the decomposition 
$$Q \simeq \bigoplus_{w \in W}\text{Pr}^{\oplus m_w}_w, \text{  }m_w \ge 0.$$
Let $y \in W.$ We have $\Hom(\text{Pr}_w, L_y) = 0$ if $w \ne y.$ Since $L_y$ is a unique simple quotient of $\text{Pr}_y$, and $\text{End}(L_y) = \C,$ we have 
$$\dim \Hom(\text{Pr}_y, L_y) = 1.$$
Hence 
$$\dim \Hom(Q, L_y) = \sum m_w \cdot \Hom(\text{Pr}_w, L_y) = m_y  $$  \Laughey

 \subsection{Recontextualizing Soergel modules from Bott-Samelson modules}
 We now take an arbitrary Coxeter system $(W,S)$ and a $W$-representation $\h$ over some field $\mathbf{k}.$ Let $R = \text{Sym}(\h^*)$ with $\deg(\h^*) = 2.$ Consider the subcategories $\overline{\mathbb{S}\text{Bim}}(\h, W)$ of Soergel bimodules and $\overline{\mathbb{BS}\text{Bim}}(\h, W)$ of Bott-Samelson bimodules. Recall that Soergel modules are obtained via these respective bimodules. We have to be a little tricky here and view $\mathbf{k}$ as a graded $R$-module with all positive degree polynomials acting by zero. 
 \subsubsection{Definition}
 Let $\underline{w} = (s_1, \ldots, s_m)$ be a reduced expression for $w.$ The \textit{Bott-Samelson module} $\overline{BS}(\underline{w})$ is the graded right $R$-module 
 $$\overline{BS}(\underline{w}) = \mathbf{k} \otimes_R BS(\underline{w}) = \mathbf{k} \otimes_R B_{s_1} \otimes_R \cdots \otimes_R B_{s_m}$$
 $$\simeq \mathbf{k} \otimes_R R \otimes_{R^{s_1}} R \otimes_{R^{s_2}} \cdots \otimes_{R^{s_m}} R(m).$$
 Now we can see more clearly how we get our Soergel modules. Explicitly, they're any graded right $R$-modules isomorphic to a finite direct sum of shifts of direct summands of Bott-Samelson modules. Our subcategory $\overline{\mathbb{S}\text{Bim}}(\h, W)$ from earlier is in fact a full subcategory of gmod-$R$. We can now consider the functor 
 $$\mathbf{k} \otimes_R (-): R\text{-gmod-}R \rightarrow \text{gmod-}R,$$
 which is an additive functor sending $\overline{BS}(\underline{w}) \mapsto \overline{BS}(\underline{w}).$ Thus it restricts to 
 $$\mathbf{k} \otimes_R (-): \mathbb{S}\text{Bim}(\h, W) \longrightarrow \overline{\mathbb{S}\text{Bim}}(\h, W).$$
 This gives us the idea that Soergel modules can equivalently be thought of as graded $R$-modules that come from direct summands of $\mathbf{k} \otimes_R B,$ where $B$ is a Soergel bimodule. To make use of this functor, we give the following

 \subsubsection{Proposition}
 If $W$ is finite, and $\h$ is reflection faithful, the assignment 
 $$\mathbf{k} \otimes_R \Hom^\bullet_{R\text{-gmod-}R}(B, B') \longrightarrow \Hom^\bullet_{\text{gmod-}R}(\mathbf{k} \otimes_R B, \mathbf{k} \otimes_R B')$$
 is an isomorphism for any Soergel bimodules $B, B'.$

 \subsubsection*{}
 Before we move into the final part of the story, recall that Soergel bimodules give the categorification of the Hecke algebra $\mathcal{H},$ and that the split Grothendieck group is a $\Z[v, v^{-1}]$-algebra with the isomorphism 
 $$[\mathbb{S}\text{Bim}(\h, W)]_\oplus \simeq \mathcal{H}.$$
 When we descend to Soergel modules, we notice that the category $\overline{\mathbb{S}\text{Bim}}(\h, W)$ is not monoidal. But we still take $\overline{\mathbb{S}\text{Bim}}(\h, W)$ to be a right module category of $\mathbb{S}\text{Bim}(\h, W),$ furnished by 
 $$(-) \otimes_R (-): \overline{\mathbb{S}\text{Bim}}(\h, W) \times \mathbb{S}\text{Bim}(\h, W) \longrightarrow \overline{\mathbb{S}\text{Bim}}(\h, W). $$
 Thus we take the split Grothendieck group $[\overline{\mathbb{S}\text{Bim}}(\h, W)]_\oplus$ to be a right $[\mathbb{S}\text{Bim}(\h, W)]_\oplus$-module. 
 
 \subsubsection*{}
 Now we can introduce the climax of our tale: proof that Soergel's conjecture implies the Kazhdan-Lusztig conjecture. Following from [Elias et. al],

 \subsection{Proving the Kazhdan-Lusztig conjecture}

\subsubsection{Restating the conjecture}
The Kazhdan-Lusztig conjecture states that for a complex semisimple Lie algebra $\g$, and any $x,y$ in the Weyl group $W,$
$$[M_y : L_x] = h_{y,x}(1)$$
where $h_{y,x}(1)$ is the value of the Kazhdan-Lusztig polynomial $h_{y,x}(v)$ at $v = 1$. By Bernstein-Gelfand-Gelfand reciprocity\ref{app5}, the Kazhdan-Lusztig conjecture is equivalent to stating 
$$(\text{Pr}_x : M_y) = h_{y,x}(1).$$
Using the correspondence $[\mathcal{O}_0] \simeq \Z[W]$ in the Grothendieck group, this means that we equivalently have 
$$[\text{Pr}_x] = b_x|_{v=1},$$
where $b_x|_{v=1}$ is the specialization at $v = 1$ of the Kazhdan-Lusztig basis element $b_x \in \mathcal{H}.$ Recall finally that Soergel's conjecture is true if for every $x \in W,$
$$\chi(B_x) = b_x.$$
We now move to our final proposition.

 \subsubsection{Proposition}\ref{ref2}
 \textit{Let $\g$ be a complex semisimple Lie algebra. Then Soergel's conjecture for the geometric representation of its Weyl group $W$ implies the Kazhdan-Lusztig conjecture for $\g.$}\\\\
 \textit{Proof:} We will begin by induction on Bruhat order. We have 
 $$\text{Pr}_{\text{id}} = M_{\text{id}},$$
 which implies that $[\text{Pr}_{\text{id}}] = \text{id} = b_{\text{id}}|_{v=1}.$ Let $x > \text{id}$ in Bruhat order. Choose $s \in S$ such that $xs < s,$ and take $x = ws,$ where $w < x.$ Then 
 $$\Theta_s\text{Pr}_w \simeq \text{Pr}_x \oplus \bigoplus_{z<x}\text{Pr}^{\oplus m_z}_z$$
 for multiplicities $m_z \ge 0.$ In the Grothendieck group, we obtain by induction 
 $$[\text{Pr}_x] = [\Theta_s\text{Pr}_w] - \sum_{z<x}m_z [\text{Pr}_z] = [\text{Pr}_w](1+s) - \sum_{z<x}m_z [\text{Pr}_z] $$
 $$= b_wb_s - \sum m_zb_z|_{v=1}.$$
 Applying $\mathbb{V}$ to $\Theta_s\text{Pr}_w,$ we get 
 $$\overline{B_w}B_s \simeq \overline{B_x} \oplus \bigoplus_{z<x}\overline{B_z}^{\oplus m_z}, $$
 and since Soergel modules decompose the same as Soergel bimodules, 
 $$B_wB_s \simeq B_x \oplus \bigoplus_{z<x}B^{\oplus m_z}_z.$$
 Taking characters and invoking Soergel's conjecture, we get 
 $$b_x = b_wb_s - \sum_{z<x}m_zb_z. $$ 
 Then by the equality 
 $$[\text{Pr}_w](1+s) - \sum_{z<x}m_z [\text{Pr}_z] =  b_wb_s - \sum m_zb_z|_{v=1}, $$ 
 we get $[\text{Pr}_x] = b_x|_{v=1},$ which is the equivalent form of the Kazhdan-Lusztig conjecture.  \Laughey

\newpage
\begin{center}
    \section{References}\label{references}
\end{center}
\begin{enumerate}

    \item\label{ref1} A. Doña Mateo, Soergel bimodules and HOMFLY-PT homology.
    University of Edinburgh https://www.maths.ed.ac.uk/~adona/files/part-iii-essay.pdf

    \item\label{ref2} B. Elias, S. Makisumi, U. Thiel, G. Williamson, Introduction to Soergel Bimodules. Springer (2021)

    \item P. Etingof, S. Gelaki, D. Nikshych, V. Ostrik, Tensor Categories. Mathematical Surveys and Monographs, vol. 205 (American Mathematical Society, Providence, RI, 2015),

    \item W. Fulton, J. Harris, Representation Theory: A First Course, Readings in Mathematics.
Graduate Texts in Mathematics, vol. 129 (Springer, New York, 1991) \\
https://doi.org/10.1007/978-1-4612-0979-9

    \item  J.E. Humphreys, Reflection Groups and Coxeter Groups. Cambridge Studies in Advanced Mathematics, vol. 29 (Cambridge University Press, Cambridge, 1990) \\ https://
doi.org/10.1017/CBO9780511623646

    \item D. Kazhdan, G. Lusztig, Representations of Coxeter groups and Hecke algebras. Invent. Math. 53(2), 165–184 (1979). https://doi.org/10.1007/BF01390031

    \item A. Kirillov Jr, Introduction to Lie groups and Lie algebras. Cambridge University Press (2010) https://doi.org/10.1017/CBO9780511755156 

    \item\label{ref9} N. Liebedinsky, A Gentle Introduction to Soergel Bimodules 1: The Basics. \\
    (2017) https://arxiv.org/abs/1702.00039

    \item E. Marberg, Positivity conjectures for Kazhdan-Lusztig theory on twisted involutions: the finite case. Journal of Algebra (2014) https://arxiv.org/abs/1306.2980
    
    \item\label{ref11} P. Polo, Construction of arbitrary Kazhdan-Lusztig polynomials in symmetric groups.
Represent. Theory 3, 90–104 (1999). https://doi.org/10.1090/S1088-4165-99-00074-6

    \item JP. Serre, Groupes de Coxeter finis : involutions et cubes. https://arxiv.org/pdf/2012.03689

   \item W. Soergel. Kategorie O, perverse Garben und Moduln uber den Koinvarianten zur Weylgruppe. J.
Amer. Math. Soc. 3 (1990)

    \item W. Soergel, Kazhdan-Lusztig-Polynome und unzerlegbare Bimoduln uber Polynomringen, J. Inst. Math.
Jussieu, 6(3) (2007)

\end{enumerate}

\newpage
\begin{center}
\section{Appendic Miscellany}\label{appendix}

\end{center}
\subsection{Computation of}\label{app1} $s^2_\alpha = e. $ 
We have
$$s_\alpha(s_\alpha(\lambda)) = s_\alpha\Big(\lambda - 2\dfrac{(\lambda, \alpha)}{(\alpha, \alpha)}\alpha\Big) = \Big(\lambda - 2\dfrac{(\lambda, \alpha)}{(\alpha, \alpha)}\alpha\Big) - \Bigg(2\dfrac{((\lambda - 2\frac{(\lambda, \alpha)}{(\alpha, \alpha)}\alpha), \alpha)}{(\alpha, \alpha)}\alpha\Bigg) $$
$$= \Big(\lambda - 2\dfrac{(\lambda, \alpha)}{(\alpha, \alpha)}\alpha\Big) -  2\dfrac{-(\lambda, \alpha)}{(\alpha, \alpha)}\alpha = \lambda - 2\dfrac{(\lambda, \alpha)}{(\alpha, \alpha)}\alpha + 2\dfrac{(\lambda, \alpha)}{(\alpha, \alpha)}\alpha$$
$$= \lambda. $$

\subsection{Geometric representation of a Coxeter system}\label{app2}
The \textit{geometric representation} of a Coxeter system $(W,S)$ is a representation $V$ of $W$ wherein $V$ is a real vector space equipped with basis $\{\alpha_s \st s \in S\}$, and bilinear form 
$$(\text{ } , \text{ }) := (\alpha_s, \alpha_t) = -\cos\dfrac{\pi}{m_{st}}$$
We are also given the familiar action of $W$ on $V$ with 
$$s(\lambda) = \lambda - 2(\lambda, \alpha_s)\alpha_s.$$

Imporant thing: \textit{$V$ is always faithful.}

\subsection{A fun but ancillary category }\label{app3}
Let $R$-gbim$_{qc}$ denote the category of graded $R$-bimodules. We equip this with a shift functor $(n)$ which sends, for each $n \in \Z,$ $M \mapsto M(n),$ where $M(n)$ denotes the usual grading on $R$-modules. Also, we give $R$-gbim$_{qc}$ the tensor product $\otimes_R$, giving $R$-gbim$_{qc}$ the structure of a monoidal category. Here the subscript "qc" stands for quasi-coherent, in that there is no assumption of finite generation of objects in the category. The tensor product is $\Z$-graded in that $(M \otimes_R N)^k $ is really the image in $M \otimes_R N$ of $\bigoplus_{i+j=k}M^i \otimes_\Z N^j.$ The monoidal identity here is the regular bimodule $R,$ and we perceive $R$-gbim$_{qc}$ as a graded category with morphisms entailed by graded $R$-bimodule maps which are homogenous of degree 0. Further information on $R$-gbim$_{qc}$ can be found in Elias et al, []

\newpage
\subsection{In aid of lemma 4.2.3,}\label{app4}
we prove the result of the following 
\subsubsection{Proposition}
\textit{Let $\mathcal{A}, \mathcal{B}$ be abelian categories, and let $F: \mathcal{A} \rightarrow \mathcal{B}, G: \mathcal{B} \rightarrow \mathcal{A} $ be functors. If $F$ is left adjoint to $G$ and $G$ is exact, then $F$ sends projective objects to projective objects.}
\\\\ \textit{Proof:}
Recall that an object $P$ in an abelian category is projective if and only if the functor $\Hom_\mathcal{A}(P, -): \mathcal{A} \rightarrow \text{Ab}$ is exact. With the above setup, let $P \in \mathcal{A}$ be a projective object. Then, for every short exact sequence 
$$\begin{tikzcd}
    0 \arrow[r] & X' \arrow[r] & X \arrow[r] & X'' \arrow[r] & 0 
\end{tikzcd}$$
the following sequence is exact: 
$$\begin{tikzcd}
    0 \arrow[r] & \Hom_\mathcal{A}(P, X') \arrow[r] & \Hom_\mathcal{A}(P, X) \arrow[r] & \Hom_\mathcal{A}(P, X'') \arrow[r] & 0. 
\end{tikzcd}$$
Since $F$ is left adjoint to $G,$ we have the natural isomorphism
$$\Hom_\mathcal{B}(F(A), B) \simeq \Hom_\mathcal{A}(A, G(B))$$
for every $A \in \mathcal{A}, B \in \mathcal{B}.$ We want to find out whether the functor $\Hom_\mathcal{B}(F(P), -): \mathcal{B} \rightarrow \text{Ab}$ is exact. By adjointness, for every object $Y$ in $\mathcal{B},$ we have the isomorphism
$$\Hom_\mathcal{B}(F(P), Y) \simeq \Hom_\mathcal{A}(P, G(Y)).$$
It follows that $\Hom_\mathcal{B}(F(P), -) \simeq \Hom_\mathcal{A}(P, G(-))$ as functors $\mathcal{B} \rightarrow \text{Ab}.$ Since $P$ is projective, $\Hom_\mathcal{A}(P, -)$ is exact as a functor $\mathcal{A} \rightarrow \text{Ab}.$ Moreover, $G$ is exact by assumption. Notice that the composition of exact functors is exact. Then, if 
$$\begin{tikzcd}
    0 \arrow[r] & Y' \arrow[r] & Y \arrow[r] & Y'' \arrow[r] & 0 
\end{tikzcd}$$
is a short exact sequence in $\mathcal{B},$ applying $G$ gives us the short exact sequence 
$$\begin{tikzcd}
    0 \arrow[r] & G(Y') \arrow[r] & G(Y) \arrow[r] & G(Y'') \arrow[r] & 0 
\end{tikzcd}$$
in $\mathcal{A}.$ Composing with $\Hom_\mathcal{A}(P,-)$ gives us the short exact sequence 
$$\begin{tikzcd}
    0 \arrow[r] & \Hom_\mathcal{A}(P, G(Y')) \arrow[r] & \Hom_\mathcal{A}(P, G(Y)) \arrow[r] & \Hom_\mathcal{A}(P, G(Y'')) \arrow[r] & 0. 
\end{tikzcd}$$ in Ab. By the isomorphism $\Hom_\mathcal{B}(F(P), -) \simeq \Hom_\mathcal{A}(P, G(-))$, it follows that 
$$\begin{tikzcd}
    0 \arrow[r] & \Hom_\mathcal{A}(F(P), Y') \arrow[r] & \Hom_\mathcal{A}(F(P), Y) \arrow[r] & \Hom_\mathcal{A}(F(P), Y'') \arrow[r] & 0 
\end{tikzcd}$$ is a short exact sequence, and hence $F(P)$ is projective in $\mathcal{B}.$ \Laughey

\subsection{Berstein-Gelfand-Gelfand Reciprocity}\label{app5}
For a complex semisimple Lie algebra $\g,$ let $\mathcal{C}$ denote the subcategory of $U(\g)$-Mod consisting of modules $V$ admitting finite dimensional weight spaces. For $V \in \mathcal{C},$ a weight $\lambda \in \h^*,$ and weight space $V[\lambda],$ we take $V[\lambda]^*$ to be all $f \in V^*$ that vanish on $V[\mu]$ where $\mu \ne \lambda.$ The dual module $\mathbf{D}(V)$ is defined by 
$$\mathbf{D}(V) = \bigoplus_{\lambda \in \h^*}V[\lambda]^*.$$
Berstein-Gelfand-Gelfand (BGG for short) reciprocity states that

\subsubsection{BGG Reciprocity}
$(\text{Pr}_\lambda : M_\mu) = [\mathbf{D}(M_\lambda) : L_\lambda] = [M_\mu : L_\lambda].$

\end{document}